\documentclass[journal,twoside]{IEEEtran} 

\usepackage{cite}
\usepackage{amsmath,amssymb,amsfonts}
\usepackage{algorithm}
\usepackage{algorithmic}
\usepackage{graphicx}
\usepackage{hyperref}
\hypersetup{hidelinks=true}
\usepackage{textcomp}

\usepackage{amsthm}
\newtheorem{assumption}{Assumption}

\newtheorem{lemma}{Lemma} 
\newtheorem{theorem}{Theorem}

\usepackage{tikz}
\usetikzlibrary{shapes.geometric, positioning, arrows.meta} 

\def\BibTeX{{\rm B\kern-.05em{\sc i\kern-.025em b}\kern-.08em
    T\kern-.1667em\lower.7ex\hbox{E}\kern-.125emX}}
    
\begin{document}
%
\tikzstyle{block} = [draw, rectangle, 
    minimum height=3em, minimum width=4em]
\tikzstyle{sum} = [draw, circle, node distance=1cm]
\tikzstyle{input} = [coordinate]
\tikzstyle{output} = [coordinate]
\tikzstyle{pinstyle} = [pin edge={to-,thin,black}]

\title{Distributed Saddle-Point Dynamics in Multilayer Networks}
\author{Christian D.~Rodr\'iguez-Camargo,
       Andr\'es F.~Urquijo-Rodr\'iguez,
        and~Eduardo~Mojica-Nava
\thanks{Manuscript received MM DD, YYYY; revised MM DD, YYYY.}
\thanks{C. D. Rodr\'iguez-Camargo is with The Atomic, Molecular, Optical and Positron Physics (AMOPP) group of the Department of Physics and Astronomy from the University College London and Programa de Investigaci\'on sobre Adquisici\'on y An\'alisis de Se\~nales (PAAS-UN) from Universidad Nacional de Colombia. E-mail: christian.rodriguez-camargo.21@ucl.ac.uk}
\thanks{A. F. Urquijo-Rodr\'iguez is with Grupo de investigación en Ingeniería e Innovación de la Universidad ECCI (GIII-ECCI) both from the Universidad ECCI, and Grupo de Superconductividad y Nanotecnolog\'ia and Grupo de \'Optica e Informaci\'on Cu\'antica both from the Department of Physics of Universidad Nacional de Colombia. E-mail: afurquijor@unal.edu.co}
\thanks{E. Mojica-Nava is with the Department of Electrical and Electronics Engineering and Programa de Investigaci\'on sobre Adquisici\'on y An\'alisis de Se\~nales (PAAS-UN) from Universidad Nacional de Colombia, Bogotá, Colombia.
E-mail: eamojican@unal.edu.co}}

\markboth{IEEE Transactions on Network Science and Engineering,~Vol.~XX, No.~X, August~20XX}%
{Shell \MakeLowercase{\textit{et al.}}: Bare Demo of IEEEtran.cls for Computer Society Journals}

\IEEEtitleabstractindextext{%

\begin{abstract}
Multilayer networks provide a more advanced and comprehensive framework for modeling real-world systems compared to traditional single-layer and multiplex networks. Unlike single-layer models, multilayer networks have multiple interacting layers, each with unique topological features. In this paper, we generalize previously developed results for distributed optimization in multiplex networks to the more general case of multilayer networks by employing a tensor formalism to represent multilayer networks and their tensor-Laplacian diffusion dynamics. Although multiplex networks are a special case of multilayer networks, where each layer has the same number of replica nodes connected one-to-one, this generalized framework removes the need for replica nodes, allowing variability in both topology and number of nodes across layers. This approach provides a fully generalized structure for distributed optimization in multilayer networks and enables more complex interlayer connections. We derive the multilayer combinatorial Laplacian tensor and extend the distributed gradient descent algorithm. We provide a theoretical analysis of the convergence of algorithms. Numerical examples validate our approach, and we explore the impact of heterogeneous layer topologies and complex interlayer dynamics on consensus time, underscoring their implications for real-world multilayer systems.
\end{abstract}

\begin{IEEEkeywords}
Distributed optimization, multilayer networks, saddle-point flow, diffusion, tensor dynamics.
\end{IEEEkeywords}}

\maketitle
\IEEEdisplaynontitleabstractindextext


\section{Introduction}
\label{sec:introduction}

\IEEEPARstart{T}{he} quantitative study of networks has provided a deep understanding of the structure and dynamics of complex systems in various fields~\cite{barabasi2016network, dorogovtsev2003evolution}, such as biology~\cite{GOSAK2018118}, neuroscience~\cite{Barabasineuro1}, engineering~\cite{8353833, 6528031}, and social sciences~\cite{Borgatti2009social}. However, the traditional single-layer network formalism does not capture the full complexity of real-world systems where multiple types of interaction occur simultaneously across distinct layers~\cite{Morenomulti2014}. Multilayer networks provide a more holistic framework for representing these interconnected systems, allowing researchers to explore a broader range of dynamics and interactions that are not apparent in single-layer or multiplex network models~\cite{bianconi2018multilayer}.

The multilayer network framework has been successfully applied in various fields. For instance, this structure has produced significant and robust results in ecological dynamics~\cite{pilosof2017multilayer}, brain networks~\cite{Mdomenico2017brain}, transport networks~\cite{Alessandretti23998083221108190transportm}, epidemiology~\cite{fujimoto2023beyond}, molecular biology~\cite{zhang2019multilayer}, and energy networks~\cite{KAZIM2024101407, toro2021multiplex}. Furthermore, multilayer networks reveal dynamical behaviors similar to those of single-layer networks, such as synchronization~\cite{PhysRevLett.114.038701}, percolation~\cite{9552851}, and diffusion and spreading processes~\cite{BOCCALETTI20141}.

Given the multidimensional nature of these processes, it is worth establishing a tensorial representation for such systems. Tensors generalize scalars, vectors, and matrices~\cite{lim2013tensors} and provide a concise framework for modeling and solving multidimensional problems. Tensors are widely used in all fields, from linear algebra and physics to engineering~\cite{6784037, qi2017tensor}. In ~\cite{DeDomenico2013MatFor} it is developed a tensor-based mathematical formalism that provides a general, robust, and comprehensive framework to represent multilayer networks with all their complexities, offering insights that remain hidden when using traditional adjacency matrices. This formalism has been applied explicitly to study disease localization~\cite{PhysRevX.7.011014}, generalization of the concept of assortativity~\cite{DEARRUDA20165}, spreading~\cite{DEARRUDA20181}, and eigenvector multicentrality~\cite{doi:10.1073/pnas.1801378116}.

Despite recent advances in multilayer networks, progress in generalizing distributed optimization problems from single-layer to multilayer networks has been limited. In conventional distributed optimization, each node typically has access to its own convex function, and the overall objective is defined as the sum of these individual functions~\cite{yang2019survey}. To solve this, each agent seeks the optimal solution to the global convex function by exchanging information with neighboring nodes and performing local computations~\cite{nedic2018distributed}. Various methods have been developed to address this problem, ranging from discrete-time gradient descent approaches~\cite{nedic2018distributed} to more recent continuous-time techniques~\cite{lin2016distributed, wu2023distributed, ma2019novel}. In particular, from a continuous-time perspective, optimization algorithms have been framed as dynamical systems~\cite{colombino2019online}, with methods based on saddle-point dynamics~\cite{cherukuri2017saddle, feijer2010stability} gaining prominence. 

{\color{black}
\subsection{Related Work}
A multilayer network is formed by several interacting networks \cite{boccaletti2014structure}, \cite{bianconi2018multilayer}. The first introduction of multilayer networks is in the context of social science to characterize the different types of social relationships that exist between the nodes of a social network \cite{bianconi2013statistical}. Currently, multilayer networks are investigated in many
fields, including diffusion dynamics and information spreading \cite{gomez2013diffusion}, \cite{perc2019diffusion}, \cite{salehi2015spreading}, epidemic spreading \cite{tian2024spreading}, neuroscience \cite{vaiana2020multilayer}, molecular biology \cite{kiani2021networks}, ecology \cite{pilosof2017multilayer}, and economy \cite{wang2018social}. Multilayer networks have also been applied effectively across various engineering domains to model and analyze complex systems with multiple interacting components such as power networks, where the interdependence between power grids and communication infrastructures is effectively captured using multilayer network models \cite{toro2021multiplex}. One layer represents the electrical grid, while another represents the communication network that monitors and controls the grid's operations. This modeling enables the analysis of cascading failures [50], [51], where a disruption in the power grid can impact the communication network and vice versa, aiding in the design of more resilient interdependent infrastructures. Another successful application is urban transportation networks including roads, railways, buses, and metro systems, which can be represented as multilayer networks where each mode of transport forms a layer \cite{du2016physics}. Information propagation is also a relevant application of multilayer networks \cite{zhuang2016information}. Although there are significant advances in other engineering applications phenomena such as percolation \cite{santoro2020optimal}, there is a lack of developments in the generalization of single-layer network distributed optimization
problems to multilayer networks.

A specific type of multilayer network is multiplex networks, where there is a one-to-one correspondence between each node across all layers. In other words, the same set of nodes is connected by different types of interactions in each layer. A multiplex network is a set of graphs composed by different layers in which the same set of nodes may be connected to each other by means of links belonging to different classes or types \cite{nicosia2013growing}. Multiplex networks are initially introduced in the context of social networks, where individuals are linked through various social relationships such as friendship, collaboration, and family, or communicate through different channels such as email, mobile phones, and instant messaging \cite{fienberg1985statistical}. Over time, the concept of multiplex networks has been extended to encompass a wide range of domains, including transportation systems \cite{sole2016congestion}, and brain networks \cite{crofts2016structure}.

Multilayer networks represent a generalization of peer-to-peer architectures, where agents are embedded within single-layer communication networks. In peer-to-peer architectures, consensus protocols have garnered significant attention and research interest due to their applications in synchronization, formation control, learning, and the coordination of multi-robot tasks, among others. These systems are inherently distributed. 

Distributed computational models and coordination mechanisms are based on local communication between agents. Such models are crucial for advancing technologies like automated transportation, smart cities, smart power grids, and robotic systems. At the heart of these technologies are devices capable of collecting data (e.g., gathering measurements from their environment), processing information, and communicating with one another to share estimates, thereby ensuring optimal and reliable system performance. Learning methods play a fundamental role in data processing, while communication between devices enables collaborative learning from the collective data within the system, without the need to share the data stored locally on the devices~\cite{9084356}.

A common framework for addressing these challenges is provided by distributed computational models, particularly those incorporating machine learning, which utilize consensus mechanisms as virtual coordinators. These mechanisms facilitate coordination and information exchange without centralized control, making them essential for scalable and efficient distributed system operation. One of the most robust tools to generate consensus mechanisms to coordinate agents are the distributed consensus-based gradient methods for convex learning problems over static graphs. The goal of these methods is to collaboratively solve the following problem:

\begin{equation}
    \min _{x \in \mathbb{R}^{p}}f(x):=\sum _{i=1}^{n}f_{i}(x)
    \label{min1}
\end{equation}
where $x$ is the global decision variable and each function $f_{i}: \mathbb{R}^{p}\rightarrow \mathbb{R}$ is convex and known by agent $i$ only. The agents are embedded in a communication network, and their goal is to obtain an optimal and consensual solution through local neighbor communications and information exchange. This local exchange is desirable in situations where the exchange of a large amount of data is prohibitively expensive due to limited communication resources. Since the problem (\ref{min1}) is over a networked system of agents $n$, the way to characterize this networks is through the definition of a graph $G$ as a tuple $G = ([m], \mathcal{E})$ of a set of agents $[n] = \{ 1, 2, 3, ..., n \}$ and a set of edges $\mathcal{E}=\{ (i,j) \quad \forall \, i,j \in [n] \}$. We write $(i,j)$ to represent an edge connecting agents $i$ and $j$. The agents want to solve the problem (\ref{min1}) collaboratively in the sense that they are willing to share some estimates with their immediate neighbors in the graph, but they are not willing (or not allowed) to share their data, which essentially means they do not reveal their functions $f_{i}$~\cite{9084356}. 

Given the communication graph $G = ([m], \mathcal{E})$, the problem is usually reformulated as follows:
\begin{eqnarray}
    \min _{x \in \mathbb{R}^{p},i\in [n]}\sum _{i=1}^{n}f_{i}(x_{i})
    \nonumber \\
    \textrm{s. t. }\, x_{i} = x_{j}, \, \forall (i,j) \in \mathcal{E} .
    \label{min2}
\end{eqnarray}

This reformulation is obtained from (\ref{min1}) by assigning a copy $x_{i}$ of the decision variable $x$ to each agent and imposing the requirement that all copies be the same. This reformulation serves as a departure point for discussing distributed methods that use consensus-type updates to enforce the equal-decision constraints in (\ref{min2})~\cite{9084356}.

Since each $f_{i}$ depends on its own variable $x_{i}$, the objective function in (\ref{min2}) is decoupled, but these variables are coupled through edge-based constraints. The idea is to distribute the problem among the agents, by allowing each agent to know its neighbors in the graph. Using the local agent knowledge of the graph and functions, every agent $i$ can solve its own local part of the overall problem. However, for the agents to collectively solve the overall problem, each agent needs to align its variables with the variables $x_{j}$ corresponding to its neighbors $j\in \mathcal{N}_{i}$, being $\mathcal{N}_{i} = \{ j : (i,j)\in \mathcal{E}\}$. As a mechanism for alignment of the variables, a consensus algorithm is used. The consensus algorithm is a distributed method that the agents can use to asymptotically agree on a decision vector. The distributed optimization algorithm updates by using a gradient of its private function, and is local in the sense that it relies on local information exchange. This method if often referred as \emph{consensus-based gradient method}, due to its use of a mixing step that resembles the distributed consensus process. There are extensions and variations of this method for systems embedded in single-layer networks. For instance, in~\cite{pu2020push} it is proposed distributed gradient-based method (push–pull and G-push–pull) where each node maintains estimates of the optimal decision variable and the average gradient of the agents’ objective functions. In~\cite{8882272} the methods are extended to time-varying graphs. There are stochastic approaches to improve the methods~\cite{9964422, 10128713, 9084360, 9131889}. Furthermore, there are methods where the focus is in the gradient tracking, that is, each agent uses an auxiliary variable to asymptotically track the gradient of the global cost in expectation~\cite{xin2019distributed, koloskova2021improved, li2020s, pu2020robust, zhou2024gradient, nguyen2023distributed, pu2021distributed}. For the interested reader, a complete and concise tutorial can be found in~\cite{9084356} and references therein.

The natural extension of the discrete-time algorithms consists of a representation to continuous-time dynamics. To construct this representation, we start with the traditional consensus equation for single-layer networks of multiagent systems using the Laplacian matrix $L$. The Laplacian matrix also called the graph Laplacian is the conventional form to represent a graph and is defined as $L = D - A$, where $D$ is the degree matrix given by $d_{ij} = \text{deg}(n_{i})\delta _{ij}$, being $\text{deg}(n_{i})$ the degree of a node $n_{i}$ that counts the number of times an edge terminates at that node, and $A$ is the adjacency matrix given by $a_{ij}=1$ if $(i,j)\in \mathcal{E}$. The understanding of these matrices will be fundamental to the generalization to tensor which we shall explain in the next sections.

Now, the single-layer consensus-based dynamics can be presented in compact form as a diffusion equation

\begin{equation} \label{consen}
\dot{x}=-Lx,    
\end{equation}
where $L$ is the graph Laplacian.  Using the consensus equation \eqref{consen}, several continuous-time saddle point dynamics have proposed to solve a distributed optimization problem in multiagent systems \cite{cherukuri2017saddle, feijer2010stability, Elia2011control, lin2016distributed, dorfler2017distributed, hauswirth2024optimization, colombino2019online}

\begin{equation}\label{opprob1}
\min _{x \in \mathbb{R}^{n}}\tilde{f}(x)=\sum _{i=1}^{n}f_{i}(x) \quad \text{s.t.} \quad
Lx = 0,
\end{equation}
where $f(x)$ is the global objective function to be minimized subject to the communication constraints between agents represented by the graph Laplacian. The corresponding saddle-point dynamics for optimization problem \eqref{opprob1}

\begin{eqnarray}
\dot{x} &=& -L x - \frac{\partial \tilde{f}(x)}{\partial x}-L\lambda,\\
\dot{\lambda}&=& Lx,
\label{lapdyn1}
\end{eqnarray}
where $\lambda$ is the vector of Lagrange multipliers \cite{bertsekas2009convex}. The main idea is to take advantage of the continuous-time representations, and the diffusion dynamics in higher-order networked systems driven by Lapacian-like operators to generalize the dynamics given by (\ref{lapdyn1}), which are useful for single-layer networks, to multilayer networks.

\subsection{Summary of Contributions}
In a previous work~\cite{rodriguezcamargo2023consensus, RODRIGUEZCAMARGO20231217}, we developed distributed optimization algorithms customized to multiplex networks, leveraging the supra-Laplacian matrix, $\mathcal{L}$, to handle diffusion dynamics within and between layers. Although these algorithms demonstrated success in modeling systems with distinct but fixed interlayer connections, such as in infrastructure networks, they are limited in their applicability to more complex multilayer systems. Specifically, multiplex networks assume identical sets of nodes across layers and restrict interlayer connections to counterpart nodes, a limitation when modeling systems with heterogeneous interlayer interactions. See Fig.~\ref{fig2}.

\begin{figure}[t]
\begin{center}
\includegraphics[width=.5\textwidth]{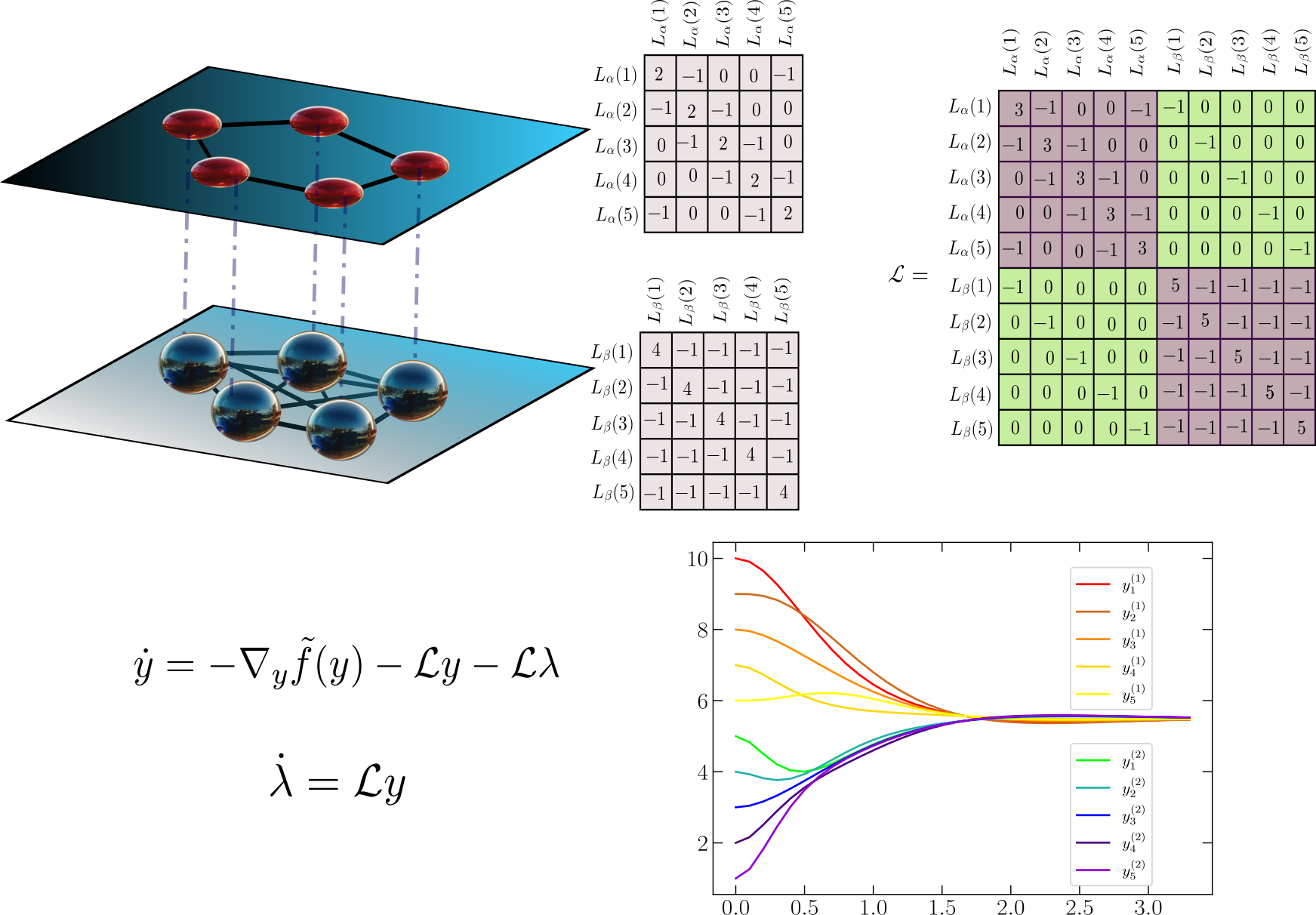}
\end{center}
\caption{Example of a multiplex network with two layers. Each layer has replica nodes, but a different topology. The connections between layers are one-by-one. On the side, we depict the Laplacians, $L_{1}$, $L_{2}$, the supra-Laplacian $\mathcal{L}$, and the optimal consensus dynamics for this system. The vector $y$ contains the state of each node within the two layers, while $\lambda$ is the vector of Lagrange multipliers. Further developments can be found in~\cite{rodriguezcamargo2023consensus, RODRIGUEZCAMARGO20231217}.}
\label{fig2}
\end{figure}

To overcome these limitations, a more generalized framework is necessary to handle complex interactions in multilayer networks. The tensorial formalism introduced in \cite{DeDomenico2013MatFor} offers a powerful approach to capture intra- and interlayer dynamics within these systems. By generalizing the network descriptors and dynamics to tensorial representations, this framework enables us to model a richer class of multilayer networks where interactions across layers are not restricted to node counterparts, but can vary in structure and strength.

\begin{figure}[th]
\begin{center}
\includegraphics[width=.36\textwidth]{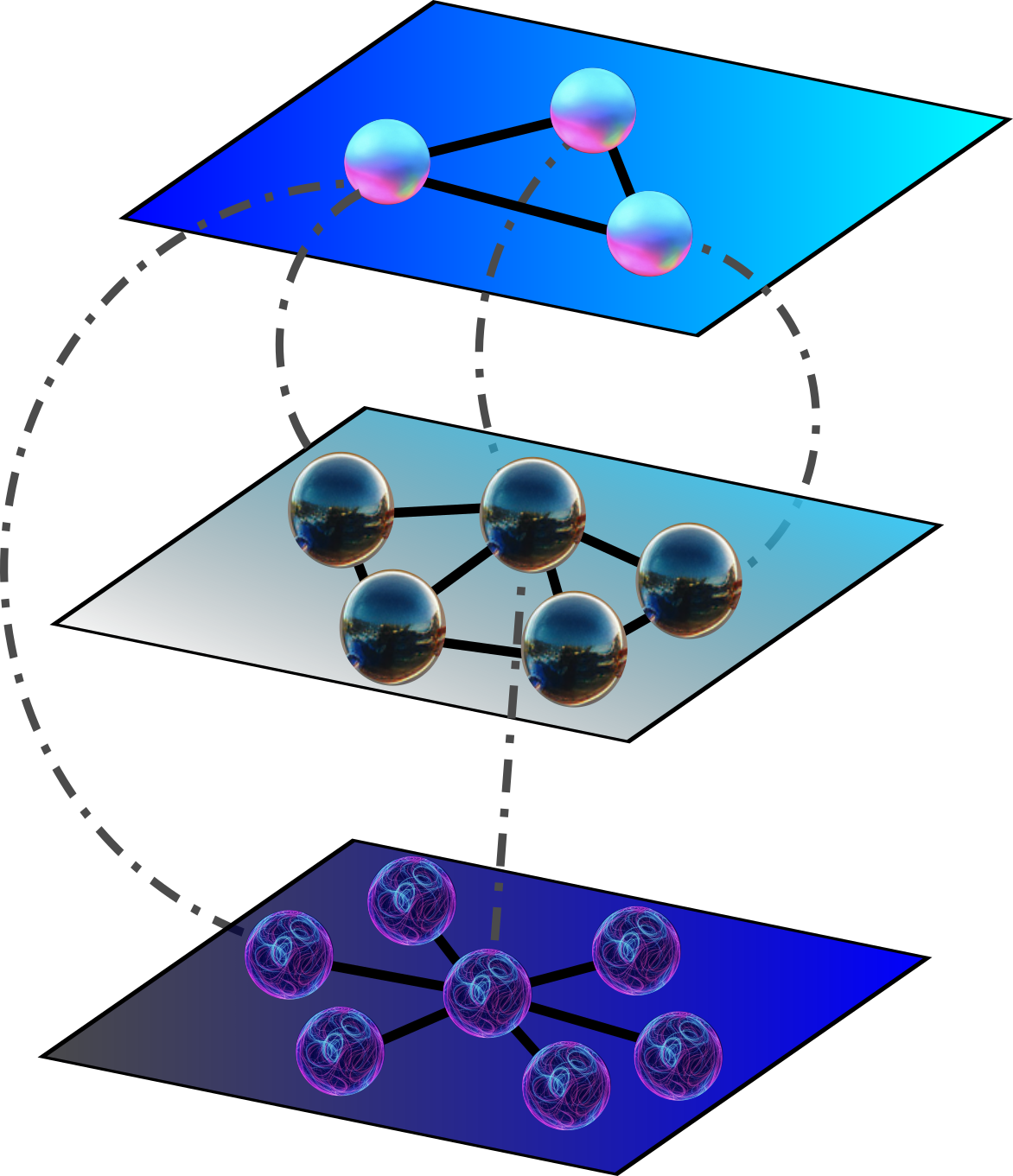}
\end{center}
\caption{Example of a multilayer network with three layers. Each layer has different number of nodes and different topology. The connections between layers may not be one-by-one.}
\label{fig4}
\end{figure}

In this paper, we extend our previous work on distributed optimization in multiplex networks by incorporating the tensorial framework for multilayer networks. This generalization allows us to handle a broader class of systems with complex interlayer connections (See Fig.~\ref{fig4}). We propose novel distributed optimization algorithms based on tensorial dynamics, analyze their convergence properties, and validate them through numerical examples. 

The main contributions of this work are as follows: (1) extending the consensus and diffusion equations to a tensor-based representation for multilayer networks, thereby generalizing the relationship between diffusion processes in multiplex and multilayer networks and consensus-driven distributed optimization algorithms; (2) formulating a distributed saddle-point algorithm tailored for convex optimization in general multilayer networks; and (3) proving the convergence of each node in every layer to the optimal value of the convex objective function. In addition, we examine the role of interlayer diffusion and its effect on optimization efficiency in multilayer networks, illustrating the applicability of our framework in real-world scenarios such as multi-energy systems. We further examine the occurrence of critical phenomena in consensus time by varying the diffusion constants. In addition, several numerical examples are presented to validate the proposed algorithms and investigate the influence of interlayer diffusion on the consensus time. We can summarize the impacts and contributions of this work in Table~\ref{tab:comparison}.

{\color{red}

\begin{table*}
\centering
\caption{Comparison of proposed approach and related references}
\label{tab:comparison}
\begin{tabular}{|l|l|l|l|l|}
\hline
\textbf{Ref.} & \textbf{Network Type} & \textbf{Optimization Method} & \textbf{Limitations} & \textbf{Our Contribution} \\ \hline
\cite{9084356, dorfler2017distributed} & \parbox{1.8cm}{Single-layer.}
 & 
\parbox{3.4cm}{Distributed Gradient descent: \\ 

$\dot{x} = -L x - \frac{\partial \tilde{f}(x)}{\partial x}-L\lambda,$\\
$\quad \dot{\lambda}= Lx$\\
} 
& Limited to single-layer networks & Extends to multilayer networks \\ \hline

\cite{rodriguezcamargo2023consensus, RODRIGUEZCAMARGO20231217} & \parbox{1.8cm}{Multiplex.\\
Fig.~\ref{fig2}} & \parbox{3.4cm}{Saddle-point dynamics: \\ 

$\dot{y}=-\nabla_{y} \tilde{f}(y)-\mathcal{L}y-\mathcal{L}\lambda ,$\\
$\dot{\lambda}= \mathcal{L}y$\\
}  & Replica nodes required  & Removes need for replica nodes \\ \hline

\cite{DeDomenico2013MatFor, DEARRUDA20181} & \parbox{1.8cm}{Multilayer.\\
Fig.~\ref{fig4}} & \parbox{3.4cm}{Tensor-based diffusion : \\ 

$\frac{dX_{\beta \Tilde{\nu}}(t)}{dt} = - \mathcal{L} ^{\alpha \Tilde{\mu}}_{\beta \Tilde{\nu}} X_{\alpha \Tilde{\mu}}(t)$\\

}   & Focused on diffusion, not optimization & Generalizes to distributed optimization \\ \hline
\end{tabular}
\end{table*}

}

This paper is organized as follows. In Section II, the fundamentals of multilayer networks are introduced. We introduce the tensor formalism given in \cite{DeDomenico2013MatFor}, and the tensor diffusion dynamics. In Section III, the distributed optimization problem is presented for general multilayer networks. We present the main result of this paper and a convergence analysis of our theorem related to the distributed primal-dual saddle-point algorithm for general multilayer networks in Section IV. In Section V numerical examples are given to illustrate the performance of the proposed approach. Conclusions and future work are presented in Section VI.
}
\section{Preliminaries}

\subsection{Tensor Representation of Multilayer Networks}

\subsubsection{Single-layer Networks}

Given a set of $N$ objects $n_{i}$, where $i=1, 2, ..., N$ and $N\in \mathbb{N}$. Each object can be associated with a state represented by a canonical vector in the vector space $\mathbb{R}^{N}$. Specifically, let $\mathbf{e}_{i}\equiv (0,...,0,1,0,...,0)^{\intercal}$, where $\intercal$ denotes the transpose operator. The column vector $\mathbf{e}_{i}$ corresponds to the object $n_{i}$, such that its $i$-th component is 1, while all other components are 0.

Since in a network, the objects $n_{i}$ have relationships among them, the vector space for describing these relationships is given by $\mathbb{R}^{N}\otimes \mathbb{R}^{N} = \mathbb{R}^{N\times N}$. Consequently, the second-order (i.e., rank-2) canonical tensors are specified by $\mathbf{E}_{ij}=\mathbf{e}_{i}\otimes \mathbf{e}_{j}^{\intercal}$, for all $i,j=1,2,...N$. Therefore, if the relations between the objects are described by the elements of a weight matrix $w_{ij}$, one can define its associate tensor as
\begin{equation}
    \mathbf{W} = \sum _{i,j=1}^{N}w_{ij}\mathbf{E}_{ij} = \sum _{i,j=1}^{N}w_{ij} \mathbf{e}_{i}\otimes \mathbf{e}_{j}^{\intercal} ,\quad \mathbf{W}\in \mathbb{R}^{N\times N} .
    \label{wtensor1}
\end{equation}

Following \cite{DeDomenico2013MatFor}, we adopt the covariant notation. In this framework, a row vector $\mathbf{a}\in \mathbb{R}^{N}$ is represented as a covariant vector $a_{\alpha}$ (where $\alpha = 1, 2, ..., N$), while its corresponding contravariant vector $a^{\alpha}$ (i.e., its dual vector) is expressed as a column vector in Euclidean space. Latin letters $i,j,...$ are used to denote vectors and tensors associated with objects (e.g., the $i$-th vector, the $(ij)$-th tensor, etc.), whereas Greek letters $\alpha, \beta, \gamma,...$ indicate the components of a vector or tensor. According to this notation, $e^{\alpha}(i)$ represents the $\alpha$-th component of the $i$-th contravariant canonical vector $\mathbf{e}_{i}\in \mathbb{R}^{N}$, and $e_{\alpha}(j)$ represents the $\alpha$-th component of the $j$-th covariant vector in $\mathbb{R}^{N}$.

This notation allows us to represent the tensor~(\ref{wtensor1}) as a linear combination of tensors in the canonical basis as follows,
\begin{equation}
    W^{\alpha}_{\beta}= \sum _{i,j=1}^{N}w_{ij}\, e^{\alpha}(i)e_{\beta}(j) = \sum _{i,j=1}^{N}w_{ij}\, E^{\alpha}_{\beta}(ij),
    \label{wtensor2}
\end{equation}
where $E^{\alpha}_{\beta}(ij) \in  \mathbb{R}^{N\times N}$ is the tensor in the canonical basis corresponding to the tensor product of the canonical vectors assigned to nodes $n_{i}$ and $n_{j}$ (i.e., $\mathbf{E}_{ij}$). The object (\ref{wtensor2}) is called the adjacency tensor. It is noticed that this tensor is of mixed type, i.e., 1-covariant and 1-contravariant.

\subsubsection{Multilayer Networks}

In multilayer networks, a node $n_{i}$ of a layer $\Tilde{h}$ can be connected to any other node $n_{j}$ in any other layer $\Tilde{k}$. The tilde symbol is introduced to distinguish indices that correspond to nodes from those that correspond to layers. A second-order interlayer adjacency tensor $C^{\alpha}_{\beta}(\Tilde{h} \Tilde{k})$ is introduced to encode the information about the relationships incorporating multiple layers. This tensor is constructed from the canonical vectors $e^{\Tilde{\mu}}(\Tilde{k})$, where $\Tilde{\mu}=1, 2, ..., M$ and $\Tilde{k} = 1, 2, ..., M$, this supposes that our multilayer network is composed by $M$ layers. These vectors are the canonical basis in the space $\mathbb{R}^{M}$. The tilde symbol on the Greek indices helps differentiate them from the Greek indices representing nodes. Note that $C^{\alpha}_{\beta}(\Tilde{k} \Tilde{k}) = W^{\alpha}_{\beta}(\Tilde{k})$, meaning that the interlayer adjacency tensor for a pair of layers representing the same layer $\Tilde{k}$ is equivalent to the intralayer adjacency tensor of that layer. 

The canonical basis (rank-2 tensors) of the space $\mathbb{R}^{M\times M}$ is then given by
\begin{equation}
    E ^{\Tilde{\mu}}_{\Tilde{\nu}} (\Tilde{h} \Tilde{k}) = e^{\Tilde{\mu}}(\Tilde{h})e_{\tilde{\nu}}(\Tilde{k}) .
\end{equation}

Following this basis, we obtain the multilayer adjacency tensor $M^{\alpha \Tilde{\mu}}_{\beta \Tilde{\nu}}$ as follows,

\begin{eqnarray}
    M^{\alpha \Tilde{\mu}}_{\beta \Tilde{\nu}} &=& \sum _{\Tilde{h}, \Tilde{k} = 1}^{M} C^{\alpha}_{\beta}(\Tilde{h} \Tilde{k}) E ^{\Tilde{\mu}}_{\Tilde{\nu}} (\Tilde{h} \Tilde{k})
    \nonumber \\
    &=& \sum _{\Tilde{h}, \Tilde{k} = 1}^{M} \left[ \sum _{i,j=1}^{N(M)}w_{ij}(\Tilde{h} \Tilde{k})  E^{\alpha}_{\beta}(ij) \right] E^{\Tilde{\mu}}_{\Tilde{\nu}}(\Tilde{h} \Tilde{k})
    \nonumber \\
    &=& \sum _{\Tilde{h}, \Tilde{k} = 1}^{M} \sum _{i,j=1}^{N(M)} w_{ij}(\Tilde{h} \Tilde{k}) \mathcal{E}^{\alpha \Tilde{\mu}}_{\beta \Tilde{\nu}} (ij\Tilde{h} \Tilde{k}),
\end{eqnarray}
where $w_{ij}(\Tilde{h} \Tilde{k})$ are real numbers indicating the intensity of the relationship, which could not be symmetric, between nodes $n_{i}$ in layer $\Tilde{h}$ and node $n_{j}$ in layer $\Tilde{k}$, and 
\begin{equation}
    \mathcal{E}^{\alpha \Tilde{\mu}}_{\beta \Tilde{\nu}} (ij\Tilde{h} \Tilde{k}) = e^{\alpha}(i)e_{\beta}(j)e^{\Tilde{\mu}}(\Tilde{h})e_{\Tilde{\nu}}(\Tilde{k})
\end{equation}
are the fourth-order (i.e., rank-4) tensors of the canonical basis in the space $\mathbb{R}^{N\times N\times M \times M}$. We have adopted the notation $N \equiv N(M)$ to remark that the number of nodes in each layer might be different. The multilayer adjacency tensor $M^{\alpha \Tilde{\mu}}_{\beta \Tilde{\nu}}$ is a versatile construct that can represent a wide range of complex relationships between nodes across different layers.

\subsection{Diffusion Processes on Multilayer Networks}

\subsubsection{Diffusion Processes on Single-layer Networks}

Using tensor notation, the diffusion processes on single-layer networks can be described as follows. Let $x_{\alpha}(t)$ be a state vector of nodes at time $t$. The diffusion equation is given by
\begin{equation}
    \frac{dx_{\beta}(t)}{dt}=D\left[ W^{\alpha}_{\beta}\, x_{\alpha}(t)- W^{\alpha}_{\gamma}u_{\alpha}e^{\gamma}(\beta)x_{\beta}(t) \right] ,
\end{equation}
where $D$ is a diffusion constant, and $u_{\alpha}$ is a 1-tensor whose components are all equal to 1, given by $u^{\alpha} = (1,...,1)^{\text{T}}\in \mathbb{R}^{N}$. The strength of a vector is obtained as $s_{\gamma} = W^{\alpha}_{\gamma}u_{\alpha}$, with this we have $s_{\gamma}e^{\gamma}(\beta)x_{\beta}=s_{\gamma}e^{\gamma}(\beta)\delta ^{\alpha}_{\beta}x_{\alpha}(t)$. Thus, the diffusion equation in single-layer networks yields:
\begin{equation}
    \frac{dx_{\beta}(t)}{dt}=-D \, L^{\alpha}_{\beta}\, x_{\alpha}(t) ,
    \label{diffeq1}
\end{equation}
where
\begin{equation}
    L^{\alpha}_{\beta}=W^{\eta}_{\gamma}u_{\eta}e^{\gamma}(\beta)\delta ^{\alpha}_{\beta}-W^{\alpha}_{\beta}
\end{equation}
is the combinatorial Laplacian tensor. This tensor is positive semidefinite, therefore, the solution of (\ref{diffeq1}) is
\begin{equation}
    x_{\beta}(t) = x_{\alpha}(0)\, \text{e}^{-DL^{\alpha}_{\beta} t} .
\end{equation}

Notice that the so-called strength tensor $\Delta ^{\alpha}_{\beta} = W^{\eta}_{\gamma}u_{\eta}e^{\gamma}(\beta)\delta ^{\alpha}_{\beta}$ is related with the degree matrix, and its components could be obtained directly such that we can express $\Delta ^{\alpha}_{\beta}$ as follows:
\begin{equation}
    \Delta ^{\alpha}_{\beta} = \sum _{i,j=1}^{N}d_{ij} e^{\alpha}(i)e_{\beta}(j),
    \label{ourdegreetensor}
\end{equation}
being $d_{ij}= d_{i}\delta _{ij}= \sum _{k}w_{ki}\delta _{ij}$. In Appendix \ref{AppendA}, we show the equivalence between these definitions. 
\subsubsection{Diffusion Processes on Multilayer Networks}

In this case, a diffusion equation for multilayer networks must incorporate terms that account for interlayer diffusion. Let $X_{\alpha \Tilde{\mu}}(t)$ denote the state tensor of the nodes in each layer at time $t$. Thus, a diffusion equation for a multilayer network is obtained as follows,
\begin{equation}
    \frac{dX_{\beta \Tilde{\nu}}(t)}{dt} = M^{\alpha \Tilde{\mu}}_{\beta \Tilde{\nu}} X_{\alpha \Tilde{\mu}}(t) - M^{\alpha \Tilde{\mu}}_{\rho \Tilde{\sigma}} U_{\alpha \Tilde{\mu}} E^{\rho \Tilde{\sigma}}(\beta \Tilde{\nu})X_{\beta \Tilde{\nu}}(t) ,
\end{equation}
being $U^{\alpha}_{\beta}=u_{\alpha}u^{\beta}$ the second-order tensor whose elements are all equal to 1.

Similar to the single-layer case, a multilayer combinatorial Laplacian tensor is presented as
\begin{equation}
    \mathcal{L} ^{\alpha \Tilde{\mu}}_{\beta \Tilde{\nu}} = M^{\eta \Tilde{\epsilon}}_{\rho \Tilde{\sigma}}U_{\eta \Tilde{\epsilon}}E^{\rho \Tilde{\sigma}}(\beta \Tilde{\nu})\delta ^{\alpha \Tilde{\mu}}_{\beta \Tilde{\nu}} - M^{\alpha \Tilde{\mu}}_{\beta \Tilde{\nu}},
    \label{multillaptensor}
\end{equation}
which is used to obtain the covariant diffusion equation for multilayer networks
\begin{equation}
    \frac{dX_{\beta \Tilde{\nu}}(t)}{dt} = - \mathcal{L} ^{\alpha \Tilde{\mu}}_{\beta \Tilde{\nu}} X_{\alpha \Tilde{\mu}}(t).
    \label{diffmulti1}
\end{equation}

Note that the solution of (\ref{diffmulti1}) is obtained as
\begin{equation}
    X_{\beta \Tilde{\nu}}(t) = X_{\alpha \Tilde{\mu}}(0) \text{e}^{- \mathcal{L} ^{\alpha \Tilde{\mu}}_{\beta \Tilde{\nu}} t}. 
\end{equation}

Therefore, it is possible to show that the tensor (\ref{multillaptensor}) is a positive semidefinite tensor.

\section{Problem Statement}
We consider the following distributed optimization problem subject to a generalized multilayer network interaction, which might be interpreted as an interdependent infrastructure

\begin{equation}
\min _{x \in \mathbb{R}}  f(x) = U^{\alpha \Tilde{\mu}}f_{\alpha \Tilde{\mu}}(x).
\label{opproblem1}
\end{equation}

To explore the assumptions and conditions of $f: \mathbb{R}\rightarrow \mathbb{R}$, let us express $f_{\alpha \Tilde{\mu}}(x)$ in the canonical basis as follows:
\begin{equation}
    f_{\alpha \Tilde{\mu}}(x) = \sum _{\Tilde{h}=1}^{M} \sum _{i=1}^{N}f_{i\Tilde{h}}(x)e_{\alpha}(i)e_{\Tilde{\mu}}(\Tilde{h}),
\end{equation}
where $f_{i\Tilde{h}}:\mathbb{R}\rightarrow \mathbb{R}$ is assumed to be Lipschitz differentiable convex cost function exclusive to agent $i$ in layer $\Tilde{h}$. It is important to note that the cost for each agent depends on the global variable $x$ across all layers. Note also that we have assumed that the decision variable $x$ is in $\mathbb{R}$, but it can be extended to $\mathbb{R}^m$. The interaction among the nodes and between layers of networks can be expressed as a cyberphysical multilayer network. To guarantee there exists a unique optimal solution $x^*$ to Problem \eqref{opproblem1}, we assume the following conditions on the gradients.

\begin{assumption}
      The gradients of the multilayer cost functions $f_{i\Tilde{h}}$ are Lipschitz continuous in every layer and interlayers $\Tilde{h}$
    \begin{equation*}
    ||\nabla f_{i\Tilde{h}}(x)-\nabla f_{i\Tilde{h}}(x')||_2 \leq L||x-x'||_2, \quad \forall i \in V(\mathcal{M}), \forall \Tilde{h} \in \mathcal{M}.
    \end{equation*}

    Where $\mathcal{M}=\{ 1,...,M \}$ is the set of layers of the multilayer network, and $V(\mathcal{M})=\{ 1,...,N(M) \}$ is the set of nodes in each layer $M\in \mathcal{M}$.
\end{assumption}

In the distributed multilayer optimization problem (\ref{opproblem1}), we should consider the heterogeneous relationship between nodes. If the multilayer is characterized by a multilayer adjacency tensor $M^{\alpha \Tilde{\mu}}_{\beta \Tilde{\nu}}$, in each layer, the weighted interaction might be different, i.e., $w_{ij}(\Tilde{h} \Tilde{k})\neq w_{ij}(\Tilde{k} \Tilde{h})$. The heterogeneity of each layer $\Tilde{h} \in \mathcal{M}$ is encoded in the object $w_{ij}(\Tilde{h} \Tilde{h})$. In addition, the generalized interlayer diffusion dynamics between layers $\Tilde{h} \in \mathcal{M}$ and $\Tilde{k} \in \mathcal{M}$ can be encoded completely in the object $w_{ij}(\Tilde{h} \Tilde{k})$.

Let $Y_{\alpha \Tilde{\mu}}$ denote the state tensor of the nodes in each layer given by
\begin{equation}
    Y_{\alpha \Tilde{\mu}} =  \sum _{\Tilde{h}=1}^{M} \sum _{i=1}^{N(M)} y_{i\Tilde{h}}e_{\alpha}(i)e_{\Tilde{\mu}}(\Tilde{h}),
\end{equation}
where $y_{i\Tilde{h}}\equiv y_{i\Tilde{h}}(t)$ is the state of the node $i$ in the layer $\Tilde{h}$, we shall proceed with the equivalent tensor-Laplacian-based consensus-constrained problem 

\begin{eqnarray}\label{opproblem3}
&&\min _{Y}\tilde{f}(Y)= U^{\alpha \Tilde{\mu}} f_{\alpha \Tilde{\mu}}(Y_{\alpha \Tilde{\mu}}) \nonumber \\
&&\text{s.t. }\quad \mathcal{L}^{\alpha \Tilde{\mu}}_{\beta \Tilde{\nu}}Y_{\alpha \Tilde{\mu}} =  0.
\end{eqnarray}
where the tensor-Laplacian constraint $\mathcal{L}^{\alpha \Tilde{\mu}}_{\beta \Tilde{\nu}}Y_{\alpha \Tilde{\mu}} =  0$ is a multilayer consensus constraint that enforces the local estimate of each agent to converge to the same consensus value. 

In the next section, we propose a saddle-point dynamics algorithm to solve the generalized distributed optimization problem with tensor-Laplacian-based constraints. 

\section{Distributed Saddle-Point Dynamics for Multilayer Networks}
In this section, we introduce the distributed continuous-time optimization algorithm for multilayer networks based on saddle-point dynamics of the tensor-Lagrangian with tensor-Laplacian constraints. Considering the constrained structure of the generalized multilayer network optimization problem \eqref{opproblem3}, we extend the saddle-point dynamics for continuous-time multiplex optimization to the generalized multilayer case using its distributed implementation.

For Problem \eqref{opproblem3}, we propose the following augmented tensor-Lagrangian function $\mathcal{L}_{m}: \mathbb{R}^{N\times M}\times \mathbb{R}^{N\times M}\rightarrow \mathbb{R}$ 
\begin{equation}
\mathcal{L}_{m}(Y,\Lambda) = U^{\alpha \Tilde{\mu}} f_{\alpha \Tilde{\mu}}(Y)+\Lambda ^{\beta \Tilde{\nu}}\mathcal{L}^{\alpha \Tilde{\mu}}_{\beta \Tilde{\nu}}Y_{\alpha \Tilde{\mu}} +\frac{1}{2}Y^{\beta \Tilde{\nu}}\mathcal{L}^{\alpha \Tilde{\mu}}_{\beta \Tilde{\nu}}Y_{\alpha \Tilde{\mu}},
\label{augmented1}
\end{equation}
we can obtain the optimal solution of \eqref{opproblem3} by solving its corresponding saddle-point problem 
\begin{equation*}
    \min_{Y}\max_{\Lambda}\mathcal{L}_{m}(Y,\Lambda),
\end{equation*}
and the associated saddle-point dynamics yields
\begin{equation}
\dot{Y}=-\nabla_{Y}\mathcal{L}_{m}(Y,\Lambda)=-\nabla_{Y} \tilde{f}(Y)-\mathcal{L}Y-\mathcal{L}\Lambda
\label{saddleflow1}
\end{equation}
\begin{equation}
\dot{\Lambda}=\nabla_{\Lambda}\mathcal{L}_{m}(Y,\Lambda)= \mathcal{L}Y ,
\label{saddleflow2}
\end{equation}

which leads to the following dynamics for each agent $i$ on each layer $\Tilde{h}$ (remarking that $y_{i\Tilde{h}}(t)$ is the state of the node $i$ in the layer $\Tilde{h}$ at time $t$) in terms of the tensor components $\beta$ and $\Tilde{\nu}$, respectively,
\begin{equation}
    \dot{Y}_{\beta \Tilde{\nu}}(t) = -\frac{\partial \Tilde{f}(Y)}{\partial Y^{\beta \Tilde{\nu}}} - \mathcal{L}^{\alpha \Tilde{\mu}}_{\beta \Tilde{\nu}}Y_{\alpha \Tilde{\mu}}(t)- \mathcal{L}^{\alpha \Tilde{\mu}}_{\beta \Tilde{\nu}}\Lambda_{\alpha \Tilde{\mu}}(t) .
    \label{saddleflow11}
\end{equation}

Analogously, for \eqref{saddleflow2} we have
\begin{equation}
    \dot{\Lambda}_{\beta \Tilde{\nu}}(t) = \mathcal{L}^{\alpha \Tilde{\mu}}_{\beta \Tilde{\nu}}Y_{\alpha \Tilde{\mu}}(t)
    \label{saddleflow22}
\end{equation}

The following lemma describes the optimality conditions for the tensor-Laplacian-based Lagrangian (\ref{augmented1}) and the steady-state of the multilayer saddle-point dynamics (\ref{saddleflow11})-(\ref{saddleflow22}).

\begin{lemma}\label{lem1}
Suppose we have a symmetric tensor-Laplacian $\mathcal{L}^{\alpha \Tilde{\mu}}_{\beta \Tilde{\nu}}=\mathcal{L}^{\beta \Tilde{\nu}}_{\alpha \Tilde{\mu}}$ associated with a weighted, undirected, and connected multilayer network. Consider convexity of the function $f_{i\Tilde{h}}$ for all $(i,\Tilde{h})$, with $i\in V=\{ 1,...,N(M) \}$ and $\Tilde{h}\in \mathcal{M}=\{ 1,...,M\}$ in tensor-Lagrangian (\ref{augmented1}). The following conditions follow:

\begin{itemize}
    \item[1] If we assume that $(Y^{\star},\Lambda ^{\star})\in \mathbb{R}^{N\times M}\times \mathbb{R}^{N\times M}$ is a saddle tensor-point of (\ref{augmented1}). Then, $(Y^{\star},\Lambda ^{\star}+\gamma U)$ is also a minimax tensor-point for any $\gamma \in \mathbb{R}$.
    \item[2] Assume $(Y^{\star},\Lambda ^{\star})\in \mathbb{R}^{N\times M}\times \mathbb{R}^{N\times M}$ is a saddle tensor-point of (\ref{augmented1}). Then, $Y^{\star} = x^{\star}U$ is the solution of the primal problem (\ref{opproblem1}) with $x^{\star} \in \mathbb{R}$.
    \item[3] A saddle-tensor-point $(Y^{\star},\Lambda ^{\star})$ exists for the multilayer augmented Lagrangian (\ref{augmented1}) satisfying
    \begin{equation*}
\mathcal{L}^{\alpha \Tilde{\mu}}_{\beta \Tilde{\nu}}\Lambda _{\alpha \Tilde{\mu}} ^{\star}+\frac{\partial \tilde{f}(Y^{\star})}{\partial Y^{\beta \Tilde{\nu}}} =0.\\
\end{equation*}
\end{itemize}
\end{lemma}

\begin{IEEEproof}
It follows from the Karush-Kuhn-Tucker conditions applied to the Lagrangian function (\ref{augmented1}) and the definition of a saddle-point in convex optimization \cite{bertsekas2009convex}.   
\end{IEEEproof}

The convergence analysis of the proposed saddle-tensor-point dynamics for multilayer networks (\ref{saddleflow11})-(\ref{saddleflow22}) is now presented. The key result of this analysis is summarized in the following theorem, which establishes the convergence of the proposed primal-dual saddle-tensor-point dynamics.

\begin{theorem}\label{thm2}
Suppose we have a twice continuously differentiable convex function $\tilde{f}: \mathbb{R}^{N\times M} \rightarrow \mathbb{R}$. The convergence of every trajectory $(Y(t), \Lambda (t))$ of the multilayer saddle-tensor-point dynamics (\ref{saddleflow11})-(\ref{saddleflow22}) to a unique saddle tensor-point $(Y^{\star}, \Lambda^{\star})$ is governed by the following conditions:

\begin{itemize}
    \item[1] $Y^{\star} = x^{\star}U$ is an optimizer of the primal optimization problem (\ref{opproblem1}) with $x^{\star} \in \mathbb{R}$.
    \item[2] $\Lambda^{\star} = \bar{\Lambda}+\text{average}(\Lambda _{0})U$, where $\text{average}(\Lambda _{0})$ is the average of the initial conditions for the tensor $\Lambda(0)$, and $\bar{\Lambda}\perp U$ satisfies 
\begin{equation*}
\mathcal{L}^{\alpha \Tilde{\gamma}}_{\beta \Tilde{\delta}}\bar{\Lambda} _{\alpha \Tilde{\gamma}} ^{\star}+\frac{\partial \tilde{f}(Y^{\star})}{\partial Y^{\beta \Tilde{\delta}}} =0.
\end{equation*}
\end{itemize}
\end{theorem}

\begin{IEEEproof}
See Appendix \ref{AppendB}.
\end{IEEEproof}

\section{Numerical Experiments}

To show the consensus dynamics in generalized multilayer networks, we have to solve the system given by (\ref{saddleflow11}) and (\ref{saddleflow22}). The explicit form of the differential equations is given in Appendix \ref{AppenC}. In this section, we consider two general multilayer situations: A simple two layers network with two links between layers, and a general four layers network with different connections.

\subsection{Multilayer Network with Two Layers}
\begin{figure}
\begin{center}
\includegraphics[width=.3\textwidth]{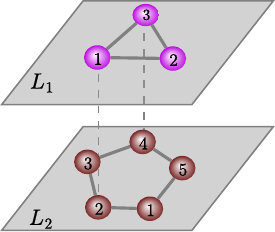}
\end{center}
\caption{Multilayer network with $M=2$ layers and $N(L_{1})=3$ nodes in layer $L_{1}$ and $N(L_{2})=5$ nodes in layer $L_{2}$. The node 1 of the layer $L_{1}$ has a connection with node 2 of layer $L_{2}$ and the node 3 of the layer $L_{1}$ has a connection with node 4 of layer $L_{2}$.}
\label{multi1}
\end{figure}

Let us start with a simple situation. We consider a multilayer network with two layers. In the first layer, we have a set of three nodes and the second layer has five nodes. In this situation, only two nodes of the first layer are connected with two nodes of the second layer. The topology is depicted in Fig.~\ref{multi1}. We study the generalized distributed optimization problem given by Problem \eqref{opproblem3} with convex function 

\begin{eqnarray}
    f_{\alpha \Tilde{\mu}}(Y_{\alpha \Tilde{\mu}}) = \frac{1}{2}Y_{\beta \Tilde{\nu}}Y^{\beta \Tilde{\nu}}\delta _{\alpha \Tilde{\mu}} + \gamma Y_{\alpha \Tilde{\mu}},
\end{eqnarray}
being $\gamma \in \mathbb{R}$ a constant which depends on each node.

\begin{figure}
\begin{center}
\includegraphics[width=.47\textwidth]{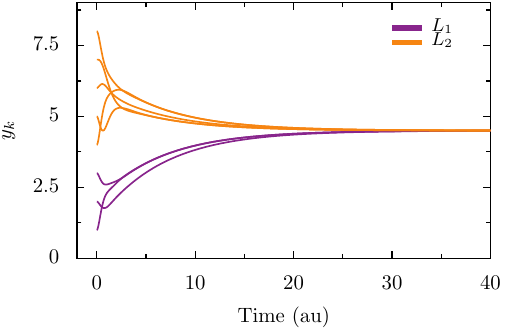}
\end{center}
\caption{Consensus dynamics for the distributed optimization problem for multilayer network with $M=2$ layers and $N(L_{1})=3$ nodes in layer $L_{1}$ and $N(L_{2})=3$ nodes in layer $L_{2}$. The interlayer diffusion constant was set to $D_{x}^{1 \rightarrow 2 } = 0.5$. }
\label{multi1r}
\end{figure}

The results of the optimal consensus are shown in Fig.~\ref{multi1r}. We can observe how the states of each node are going to the optimal solution, showing the validity of these tensor-Laplacian-driven dynamics. 

\subsection{Multilayer Network with Four Layers with Non-successive Interlayer Connections}

\begin{figure}
\begin{center}
\includegraphics[width=.5\textwidth]{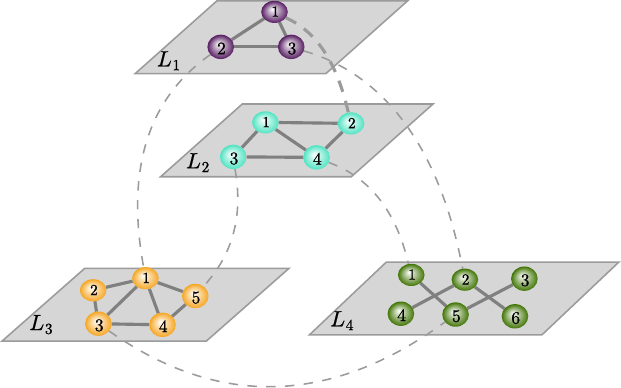}
\end{center}
\caption{Multilayer network with $M=4$ layers and $N(L_{1})=3$, $N(L_{2})=4$, $N(L_{3})=5$, $N(L_{6})=6$.}
\label{multi2}
\end{figure}

\begin{figure}
\begin{center}
\includegraphics[width=.45\textwidth]{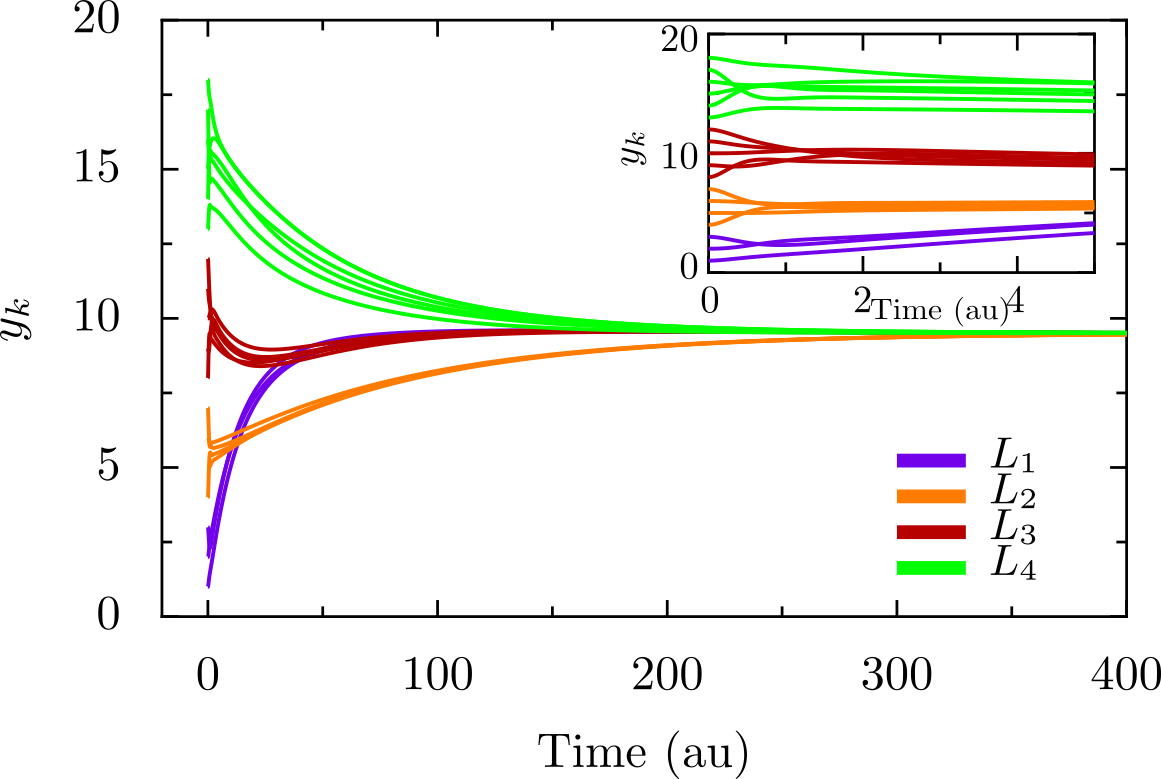}
\end{center}
\caption{Consensus dynamics for multilayer network with $M=4$ layers and $N(L_{1})=3$, $N(L_{2})=4$, $N(L_{3})=5$, $N(L_{6})=6$. The interlayer diffusion coefficientes were set to: $D_{x}^{1 \rightarrow 2} = 0.1$, $D_{x}^{1 \rightarrow 3} = 0.5$, $D_{x}^{1 \rightarrow 4} = 0.2$, $D_x^{ 2 \rightarrow 3 } = 0.1$, $D_{x}^{2 \rightarrow 4 }$, and $D_{x}^{3 \rightarrow 4} = 0.3$.}
\label{multi2r}
\end{figure}

In this case, we consider a multilayer network with four layers with different number of nodes in each layer and non-successive interlayer connections. In the first layer, we have 3 nodes. The second layer has 4 nodes. The third layer has 5 nodes. The fourth layer has 6 nodes. Each node of the first layer is connected with one node of the other three layers. Two nodes of the second layer are connected, each one, with one node of layers 3 and 4. The topology is depicted in Fig.~\ref{multi2}. The optimization problem is the same as in the previous example. The results are shown in Fig.~\ref{multi2r}. We can observe that the states of the nodes are taking more time to reach the optimal consensus in comparison with the previous situation. The inset of Fig.~\ref{multi2r} gives a zoom to early times. From this, we can evidence how the states of each layer search their intra-layer consensus, and after they do move forward to the global consensus.

\begin{figure}
\begin{center}
\includegraphics[width=.5\textwidth]{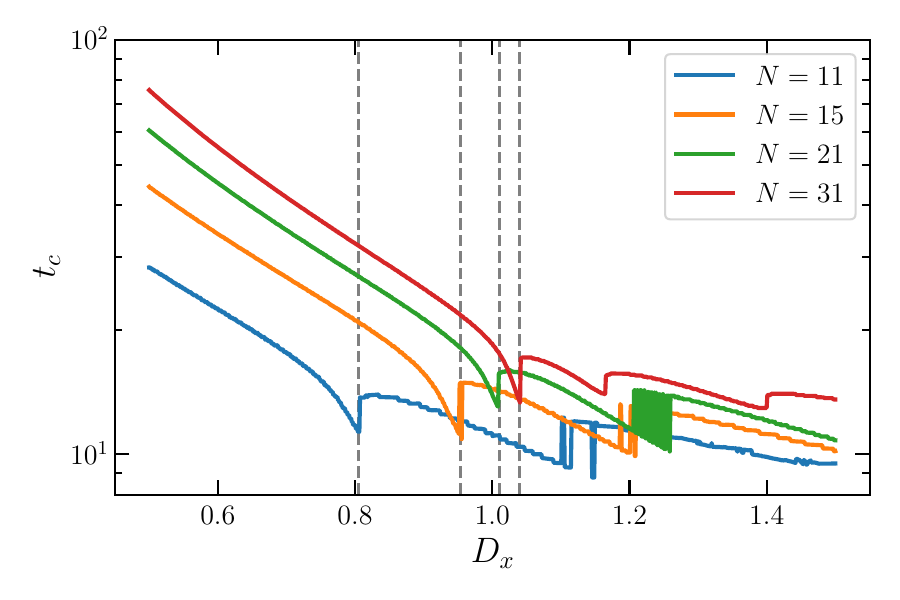}
\end{center}
\caption{Consensus time as a function of inter-layer diffusion constants for multilayer network with $M=4$ layers and $N(L_{1})=3$, $N(L_{2})=4$, $N(L_{3})=5$, $N(L_{6})=6$. }
\label{multi2rc}
\end{figure}

In Fig.~\ref{multi2rc} we plot the behavior of the consensus time $t_{C}$ as a function of the inter-layers diffusion constants. As was observed for multiplexes~\cite{PhysRevLett.110.028701, rodriguezcamargo2023consensus, RODRIGUEZCAMARGO20231217}, these dynamics are exhibiting discontinuities or criticality. That is to say, the inter-layer diffusion constants may serve as control parameters for the consensus times. 

\subsection{Multiplex Networks with Two Layers and Five Replica Nodes}

\begin{figure}
\begin{center}
\includegraphics[width=.3\textwidth]{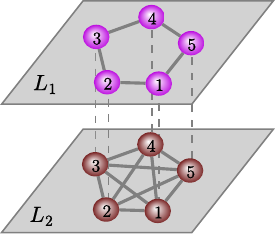}
\end{center}
\caption{Multiplex network with two layers and 5 replica nodes.}
\label{multi3}
\end{figure}

\begin{figure}
\begin{center}
\includegraphics[width=.45\textwidth]{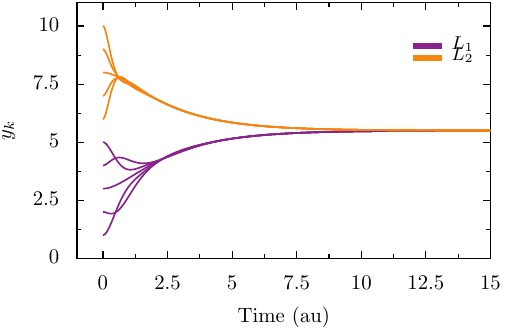}
\end{center}
\caption{Consensus dynamics for a multiplex network with two layers and 5 replica nodes. The interlayer diffusion constant was set to $D_{x}^{1 \rightarrow 2 } = 0.4$.}
\label{multi3r}
\end{figure}

For completeness, since multiplexes are a particular case of multilayer networks, we run our algorithm in a known situation where we have a multiplex of two layers and 5 replica nodes. In Fig.~\ref{multi3} we depict the topology. The results are shown in Fig.~\ref{multi3r} As expected, the results coincide with those obtained with the supra-Laplacian matrix defined in \cite{PhysRevLett.110.028701}.

\section{Conclusion}

Multilayer networks represent complex systems with multiple interacting layers, each with a unique topology and non-trivial interlayer connections. These networks have emerged as powerful tools for modeling such systems and understanding interactions between different networks. They offer a holistic framework that enables researchers to explore dynamics and interactions that are not visible in single-layer or multiplex models. Additionally, the tensor formalism provides a natural mathematical structure to capture the complexities inherent in the multidimensional nature of multilayer processes.

In this paper, we develop a distributed optimization algorithm tailored for generalized multilayer networks within a tensorial framework. By extending the connection between diffusion and consensus dynamics and leveraging control approaches used in multiplex networks, we propose a distributed tensorial primal-dual saddle-point algorithm that achieves fast convergence and robustness. Similarly to multiplex networks, we find that intra- and interlayer diffusion constants serve as control parameters for optimal consensus time. Each layer first reaches a local consensus before achieving global consensus across layers.

We also identify a critical phenomenon related to the interlayer diffusion constant, which affects consensus time. This finding extends previous research on diffusion dynamics in multiplex networks, where a phase transition was observed, associated with the second eigenvalue of the supra-Laplacian matrix. Here, phase transitions are tied to the tensor Laplacian.

With this mathematical formulation, which simplifies the general treatment of network generalization, we anticipate that the tensorial framework will enable the design of additional algorithms for real-world complex systems such as power grids, electronic circuits, brain networks, and molecular biology complexes. This generalization paves the way for exploring network properties that remain hidden within more restricted frameworks.

Possible extensions of this work involve obtaining a diffusion operator for more general higher-order networks and the application of the multilayer optimization proposed approach~\cite{battiston2020networks}. A natural continuation of this work is to study hypergraphs consensus dynamics~\cite{PhysRevE.101.062313} and consensus in simplicial complexes~\cite{majhi2022dynamics, Gong_2024}. On the other hand, topological Dirac operators offer a generalization scenario, which could be exploited to reveal new properties in complex systems, including new consensus dynamics~\cite{carletti2024globaltopologicaldiracsynchronization}. Another potential extension involves the introduction of disorder~\cite{PhysRevE.104.034102, PhysRevE.90.012146} and the application of graph zeta functions to analyze diffusion dynamics in higher-order networks~\cite{saldivar2020functional}.

\appendices
\section{Explicit Euclidean Components of $\Delta ^{\alpha}_{\beta}$ and $\Delta ^{\alpha \tilde{\mu}}_{\beta \tilde{\nu}}$ and Equivalence between Definitions} \label{AppendA}
\label{appcomponents}
Let  us start with the definition given in \cite{battiston2020networks} 
\begin{equation}
    \Delta ^{\alpha}_{\beta} = W ^{\eta}_{\gamma}u_{\eta}e^{\gamma}(\beta)\delta ^{\alpha}_{\beta} .
\end{equation}

Notice that in the indicial notation it is possible to have 
\begin{equation}
    e^{\alpha}(i) = \delta ^{\alpha} _{(i)};\qquad e_{\beta}(\alpha) = \delta _{\beta (\alpha)} .
\end{equation}

Thus,
\begin{eqnarray*}
    \Delta ^{\alpha}_{\beta} &=& W ^{\eta}_{\gamma}u_{\eta}e^{\gamma}(\beta)\delta ^{\alpha}_{\beta}\\
    \nonumber 
    &=& W ^{\eta}_{\gamma}u_{\eta}\sum _{i} \delta ^{\gamma} _{(i)}\delta ^{\alpha} _{(i)} \delta _{\beta (i)}\\
    \nonumber
    &=& \sum _{i,j} W ^{\eta}_{\gamma}u_{\eta} e ^{\gamma}(i)\delta _{ij}e ^{\alpha}(i) e _{\beta}(j)\\
    \nonumber
    &=& \sum _{i,j} k_{\gamma}e ^{\gamma}(i)\delta _{ij}e ^{\alpha}(i) e _{\beta}(j)\\
    \nonumber
    &=& \sum _{i,j} \left[ \sum _{k,l}w_{kl}\, e_{\gamma}(l) \right] e ^{\gamma}(i)\delta _{ij}e ^{\alpha}(i) e _{\beta}(j)
    \\
    \nonumber
    &=& \sum _{i,j}\sum _{k,l}w_{kl}\, \delta _{il} \delta _{ij}e ^{\alpha}(i) e _{\beta}(j)\\
    \nonumber
    &=& \sum _{i,j}\sum _{k}w_{ki}\, \delta _{ij}e ^{\alpha}(i) e _{\beta}(j)\\
    \nonumber
    &=& \sum _{i,j}d_{i}\, \delta _{ij}e ^{\alpha}(i) e _{\beta}(j)\\
    \nonumber
    &=& \sum _{i,j}d_{ij}e ^{\alpha}(i) e _{\beta}(j) .
\end{eqnarray*}

The last line is Equation (\ref{ourdegreetensor}).

\subsection{Explicit Components Definition for $\Delta ^{\alpha \tilde{\mu}}_{\beta \tilde{\nu}}$}
Using a similar procedure we can obtain a relation between the components of the multistrength tensor $\Delta ^{\alpha \tilde{\mu}}_{\beta \tilde{\nu}}$, noted as $d_{ij}(\tilde{h}\tilde{k})$, and the components of the multilayer adjacency tensor $\Delta ^{\alpha \tilde{\mu}}_{\beta \tilde{\nu}}$, given by $w_{ij}(\tilde{h}\tilde{k})$.

\begin{eqnarray*}
    \Delta ^{\alpha \tilde{\mu}}_{\beta \tilde{\nu}} &=& M^{\eta \tilde{\epsilon}}_{\rho \tilde{\sigma}}U_{\eta \tilde{\epsilon}} E^{\rho \tilde{\sigma}}(\beta \tilde{\nu})\delta ^{\alpha \tilde{\mu}}_{\beta \tilde{\nu}}\\
    \nonumber
    &=& M^{\eta \tilde{\epsilon}}_{\rho \tilde{\sigma}}u_{\eta} u_{\tilde{\epsilon}} e^{\rho}(\beta) e^{\tilde{\sigma}}(\tilde{\nu})\delta ^{\alpha \tilde{\mu}}_{\beta \tilde{\nu}}\\
    \nonumber
    &=& M^{\eta \tilde{\epsilon}}_{\rho \tilde{\sigma}}u_{\eta} u_{\tilde{\epsilon}} \delta ^{\rho}_{(\beta)} \delta^{\tilde{\sigma}}_{(\tilde{\nu})}\delta ^{\alpha}_{\beta} \delta ^{\tilde{\mu}}_{\tilde{\nu}}\\
    \nonumber
    &=&\sum _{\tilde{h},\tilde{k}}\sum _{i,j} K_{\rho \tilde{\sigma}} e^{\rho}(i)e^{\alpha}(i)e_{\beta}(j)\delta _{ij}e^{\tilde{\sigma}}(\tilde{h})e^{\tilde{\mu}}(\tilde{h})e_{\tilde{\nu}}(\tilde{k})\delta_{\tilde{h}\tilde{k}}
    \\
    \nonumber
    &=&\sum _{\tilde{h},\tilde{k}}\sum _{i,j} k_{i\tilde{h}}\delta _{ij}\delta_{\tilde{h}\tilde{k}} e^{\alpha}(i)e_{\beta}(j)e^{\tilde{\mu}}(\tilde{h})e_{\tilde{\nu}}(\tilde{k}),
\end{eqnarray*}
where the coefficients $k_{i\tilde{h}}$ are the elements of the multidegree tensor defined by
\begin{eqnarray}
    K_{\rho \tilde{\sigma}} &=& M^{\eta \tilde{\epsilon}}_{\rho \tilde{\sigma}}U_{\eta \tilde{\epsilon}}\\
    \nonumber
    &=& \sum _{\tilde{h}}\sum _{i}k_{i\tilde{h}}e_{\rho}(i)e_{\tilde{\sigma}}(\tilde{h}) ,
\end{eqnarray}
being
\begin{equation}
    k_{i\tilde{h}} = \sum _{\tilde{k}}\sum _{j}w_{ij}(\tilde{h}\tilde{k}) .
\end{equation}

Therefore, the final form of the multistrength tensor yields,

\begin{eqnarray}
    \Delta ^{\alpha \tilde{\mu}}_{\beta \tilde{\nu}} &=& \sum _{\tilde{h},\tilde{k}}\sum _{i,j}  d_{ij}(\Tilde{h} \Tilde{k}) \mathcal{E}^{\alpha \Tilde{\mu}}_{\beta \Tilde{\nu}} (ij\Tilde{h} \Tilde{k})\\
    \nonumber
    &=& \sum _{\tilde{h},\tilde{k}}\sum _{i,j} \left[ \sum _{\tilde{n}}\sum _{m} w_{im}(\tilde{h}\tilde{n}) \right] \delta _{ij}\delta_{\tilde{h}\tilde{k}}\mathcal{E}^{\alpha \Tilde{\mu}}_{\beta \Tilde{\nu}} (ij\Tilde{h} \Tilde{k}) .
\end{eqnarray}

\section{Proof of Theorem \ref{thm2}} \label{AppendB}
\begin{IEEEproof}
First, consider the convergence analysis for the tensor points $(Y^{\star}_{\beta \Tilde{\nu}}, \Lambda^{\star}_{\beta \Tilde{\nu}})$. Let $\tilde{Y}_{\beta \Tilde{\nu}}=Y_{\beta \Tilde{\nu}}-Y^{\star}_{\beta \Tilde{\nu}}$ and $\tilde{\Lambda}_{\beta \Tilde{\nu}}=\Lambda_{\beta \Tilde{\nu}} -\Lambda^{\star}_{\beta \Tilde{\nu}}$.  In the new variables, the saddle-point flow (\ref{saddleflow1})-(\ref{saddleflow2}) yields

\begin{equation}
\dot{\tilde{Y}}_{\beta \Tilde{\nu}}=-\frac{\partial \tilde{f}(Y)}{\partial Y^{\beta \Tilde{\nu}}} + \frac{\partial \tilde{f}(Y^{\star})}{\partial Y^{\beta \Tilde{\nu}}} -\mathcal{L}^{\alpha \Tilde{\mu}}_{\beta \Tilde{\nu}}\tilde{Y}_{\alpha \Tilde{\mu}}-\mathcal{L}^{\alpha \Tilde{\mu}}_{\beta \Tilde{\nu}}\tilde{\Lambda}_{\alpha \Tilde{\mu}},
\label{saddleflow3}
\end{equation} 

\begin{equation}
\dot{\tilde{\Lambda}}_{\beta \Tilde{\nu}}=\mathcal{L}^{\alpha \Tilde{\mu}}_{\beta \tilde{\nu}}\tilde{Y}_{\alpha \Tilde{\mu}} .
\label{saddleflow4}
\end{equation}

Consider the following quadratic candidate Lyapunov function
\begin{equation}
V(\tilde{Y},\tilde{\Lambda})=\frac{1}{2}\tilde{Y}^{\alpha \tilde{\mu}}\tilde{Y}_{\alpha \tilde{\mu}}+\frac{1}{2}\tilde{\Lambda}^{\alpha \tilde{\mu}}\tilde{\Lambda}_{\alpha \tilde{\mu}}.
\label{lyapunov1}
\end{equation}

Then, the temporal derivative of (\ref{lyapunov1}) along the trajectories of (\ref{saddleflow1})-(\ref{saddleflow2}) is
\begin{equation}
\dot{V}(\tilde{Y},\tilde{\Lambda})=-\tilde{Y}^{\alpha \tilde{\mu}}\frac{\partial \tilde{f}(Y)}{\partial Y ^{\alpha \tilde{\mu}}} +\tilde{Y}^{\alpha \tilde{\mu}}\frac{\partial \tilde{f}(Y^{\star})}{\partial Y ^{\alpha \tilde{\mu}}}-\tilde{Y}^{\alpha \tilde{\mu}}\mathcal{L}^{\beta \tilde{\nu}}_{\alpha \tilde{\mu}}\tilde{Y}_{\beta \tilde{\nu}}.
\end{equation}

Given the property that the gradient serves as a global under-estimator, we can conclude that
\begin{equation}
\tilde{f}(Y^{\star})\geq \tilde{f}(Y)+(Y^{\star \, \alpha \tilde{\mu}}-Y^{\alpha \tilde{\mu}})\frac{\partial \tilde{f}(Y)}{\partial Y ^{\alpha \tilde{\mu}}}.
\label{ineqproof}
\end{equation}

From inequality (\ref{ineqproof}), and if $\tilde{f}$ is strictly convex, we have
\begin{equation*}
-\tilde{Y}^{\alpha \tilde{\mu}}\frac{\partial \tilde{f}(Y)}{\partial Y ^{\alpha \tilde{\mu}}}+\tilde{Y}^{\alpha \tilde{\mu}}\frac{\partial \tilde{f}(Y^{\star})}{\partial Y ^{\alpha \tilde{\mu}}} < 0.
\end{equation*}

Now, since the Laplacian tensor is semi-positive definite, we obtain that
\begin{equation*}
\dot{V}(\tilde{Y},\tilde{\Lambda})=-\tilde{Y}^{\alpha \tilde{\mu}}\frac{\partial \tilde{f}(Y)}{\partial Y ^{\alpha \tilde{\mu}}} +\tilde{Y}^{\alpha \tilde{\mu}}\frac{\partial \tilde{f}(Y^{\star})}{\partial Y ^{\alpha \tilde{\mu}}}-\tilde{Y}^{\alpha \tilde{\mu}}\mathcal{L}^{\beta \tilde{\nu}}_{\alpha \tilde{\mu}}\tilde{Y}_{\beta \tilde{\nu}} <0.
\end{equation*}

Nevertheless, note that $\dot{V}(\tilde{Y},\tilde{\Lambda})=0$ in the set $E=\{ (\tilde{Y},\tilde{\Lambda})|\tilde{Y}=0 \}$. 

To demonstrate the global asymptotic stability of the equilibrium tensor point $(Y^{\star},\Lambda ^{\star})$, we apply LaSalle's invariance principle, showing that there are no trajectories in $E$ other than the equilibrium points $\tilde{Y}=0$ and $\tilde{\Lambda}=0$. It is observed that if $\tilde{Y}=0$ the system (\ref{saddleflow3})-(\ref{saddleflow4}) casts to
\begin{equation*}
\dot{\tilde{Y}}^{\alpha \tilde{\mu}}=-\mathcal{L}^{\alpha \tilde{\mu}}_{\beta \tilde{\nu}}\tilde{\Lambda}^{\beta \tilde{\nu}}
\end{equation*}
\begin{equation*}
\dot{\tilde{\Lambda}}^{\alpha \tilde{\mu}}=0.
\end{equation*}

Recall that the tensor Laplacian $\mathcal{L}^{\alpha \tilde{\mu}}_{\beta \tilde{\nu}}$ is semi-positive definite and by Lemma \ref{lem1} the solution is $(0, \varsigma \, U_{\alpha \tilde{\mu}})$ or $(0,0)$. According to Lemma \ref{lem1}, we can identify as $Y^{\star}_{\alpha \tilde{\mu}}=x^{\star}U_{\alpha \tilde{\mu}}$ the optimal tensor points $(Y^{\star}, \Lambda^{\star})$. Moreover, $\Lambda^{\star}_{\alpha \tilde{\mu}} = \bar{\Lambda}_{\alpha \tilde{\mu}}+\varsigma \, U_{\alpha \tilde{\mu}}$, where $\varsigma \in \mathbb{R}$ and $\bar{\Lambda}_{\alpha \tilde{\mu}}\perp U_{\alpha \tilde{\mu}}$ satisfies
\begin{equation*}
\mathcal{L}^{\alpha \tilde{\mu}}_{\beta \tilde{\nu}}\bar{\Lambda}_{\alpha \tilde{\mu}}+\frac{\partial \tilde{f}(Y^{\star})}{\partial Y ^{\alpha \tilde{\mu}}}=0.
\end{equation*}

Examining the dual dynamics (\ref{saddleflow4}), we can observe that $U^{\alpha \tilde{\mu}}\dot{\Lambda}_{\alpha \tilde{\mu}}=0$. Therefore, we obtain $\text{average}(\Lambda (t))=\text{average}(\Lambda_{0})$ $\forall t \geq 0$. As a consequence, the convergence to a single tensor saddle-point $(Y^{\star}, \Lambda ^{\star})$ of the flow (\ref{saddleflow3})-(\ref{saddleflow4}) is guaranteed, which satisfies $Y^{\star}_{\alpha \tilde{\mu}}=x^{\star}U_{\alpha \tilde{\mu}}$ and $\Lambda^{\star}_{\alpha \tilde{\mu}} = \bar{\Lambda}_{\alpha \tilde{\mu}}+ \text{average}(\Lambda_{0}) U_{\alpha \tilde{\mu}}$.
\end{IEEEproof}

\section{Explicit Form of the Generalized Consensus Dynamics}
\label{AppenC}
Since $Y_{\alpha \Tilde{\mu}}$ denotes the state tensor of the nodes in each layer given by
\begin{equation}
    Y_{\alpha \Tilde{\mu}} =  \sum _{\Tilde{h}=1}^{M} \sum _{i=1}^{N(M)} y_{i\Tilde{h}}e_{\alpha}(i)e_{\Tilde{\mu}}(\Tilde{h}),
\end{equation}
where $y_{i\Tilde{h}}\equiv y_{i\Tilde{h}}(t)$ is the state of the node $i$ in the layer $\Tilde{h}$, we define the temporal derivative of $Y_{\alpha \Tilde{\mu}}$ as
\begin{equation}
    \dot{Y}_{\alpha \Tilde{\mu}}=  \sum _{\Tilde{h}=1}^{M} \sum _{i=1}^{N(M)} \dot{y}_{i\Tilde{h}}e_{\alpha}(i)e_{\Tilde{\mu}}(\Tilde{h}) .
\end{equation}

On the other hand, the 4-rank tensor regarding the generalized gradient is defined as
\begin{equation}
    \frac{\partial f_{\alpha \tilde{\mu}}(Y)}{\partial Y^{\beta \tilde{\nu}}} = \sum _{\Tilde{h}, \Tilde{k} = 1}^{M} \sum _{i,j=1}^{N(M)} \frac{\partial f_{i\tilde{h}}}{\partial y_{j\tilde{k}}}\mathcal{E}_{\alpha \beta\Tilde{\mu}\Tilde{\nu}} (ij\Tilde{h} \Tilde{k}) .
\end{equation}

Therefore, the gradient of the contraction $\tilde{f} = U^{\alpha \Tilde{\mu}} f_{\alpha \Tilde{\mu}}$ is given by
\begin{equation}
    \frac{\partial \Tilde{f}(Y)}{\partial Y^{\beta \Tilde{\nu}}} = \sum _{\Tilde{h}, \Tilde{k} = 1}^{M} \sum _{i,j=1}^{N(M)} \frac{\partial f_{j\tilde{k}}}{\partial y_{i\tilde{h}}}e_{\beta}(i) e_{\tilde{\nu}}(\tilde{h}) .
\end{equation}

Finally, the contraction of the Laplacian tensor and any other 2-rank tensor $T_{\alpha \tilde{\mu}}$ yields
\begin{equation}
    \mathcal{L}^{\alpha \tilde{\mu}}_{\beta \tilde{\nu}}T_{\alpha \tilde{\mu}} = \sum _{\Tilde{h}, \Tilde{k} = 1}^{M} \sum _{i,j=1}^{N(M)} l_{ij}(\tilde{h}\tilde{k})t_{j\tilde{k}} e_{\beta}(i) e_{\tilde{\nu}}(\tilde{h}),
\end{equation}
where the elements $l_{ij}(\tilde{h}\tilde{k})$ are given by
\begin{eqnarray}
    l_{ij}(\tilde{h}\tilde{k}) &=& d_{ij}(\tilde{h}\tilde{k}) - w_{ij}(\tilde{h}\tilde{k}) 
    \nonumber \\
    &=& \left[ \sum _{\Tilde{n} = 1}^{M} \sum _{m=1}^{N(M)} w_{mi}(\tilde{n}\tilde{h}) \right] \delta _{ij}\delta_{\tilde{h}\tilde{k}} - w_{ij}(\tilde{h}\tilde{k}) .
    \nonumber \\
\end{eqnarray}
    
Putting all these elements together in the expression (\ref{saddleflow11}) we obtain that the system of differential equations for the evolution of the state $y_{i\tilde{h}}(t)$ of the node $i$ in the layer $\tilde{h}$ is the following,
\begin{eqnarray}
    \dot{y}_{i\Tilde{h}}(t) &=& \sum _{\Tilde{k} = 1}^{M} \sum _{j=1}^{N(M)}\frac{\partial f_{j\tilde{k}}}{\partial y_{i\tilde{h}}} 
    \nonumber \\
    &-&  \sum _{\Tilde{k} = 1}^{M} \sum _{j=1}^{N(M)} l_{ij}(\tilde{h}\tilde{k}) y_{j\tilde{k}}(t) 
    - \sum _{\Tilde{k} = 1}^{M} \sum _{j=1}^{N(M)} l_{ij}(\tilde{h}\tilde{k}) \lambda_{j\tilde{k}}(t) .
    \nonumber \\
\end{eqnarray}
 \section*{Acknowledgments}
This paper was partially supported by Minciencias Grant number CT
542-2020, \emph{Programa de Investigación en Tecnologías Emergentes para Microrredes Electricas Inteligentes con Alta Penetración de Energías Renovables} and the \emph{VIII Convocatoria para el Desarrollo y Fortalecimiento de los Grupos de Investigaci\'on en Uniminuto} with code C119-173, the \emph{Convocatoria de investigación para prototipado de tecnologías que promueven el cuidado o la restauración del medioambiente} with code CPT123-200-5220, and Industrial Engineering Program from the Corporaci\'on Universitaria Minuto de Dios (Uniminuto, Colombia). We thank the Engineering and Physical Sciences Research Council (EPSRC) (Grants No. EP/R513143/1 and No. EP/T517793/1) for financial support.




\bibliographystyle{IEEEtran}
\bibliography{Multiplexbib}

%

%

\begin{IEEEbiographynophoto}{Christian David Rodr\'iguez-Camargo}
Received a B.S. degree in physics from Universidad Nacional de Colombia in 2014, an M.Sc. degree in physics from the Brazilian Center for Research in Physics (CBPF) in 2016. He does research in relativistic quantum information, quantum biology, quantum field theory in curved spaces, smart grits and statistical field theory in network science within researching collaboration at The Atomic, Molecular, Optical and Positron Physics (AMOPP) group of the Department of Physics and Astronomy from the University College London and Programa de Investigaci\'on sobre Adquisici\'on y An\'alisis de Se\~nales (PAAS-UN) from Universidad Nacional de Colombia.
\end{IEEEbiographynophoto}

\begin{IEEEbiographynophoto}{Andr\'es F. Urquijo-Rodr\'iguez}
Received a B.S. and M.Sc. degree in physics from Universidad Nacional de Colombia, and he is currently pursuing his Ph.D. from Universidad Nacional de Colomba. He does research in network science, quantum computing, and condensed matter physics in low dimensional systems into polaritonic systems. He is a part-time professor at Pontificia Universidad Javeriana, Bogotá, Colombia.
\end{IEEEbiographynophoto}

\begin{IEEEbiographynophoto}{Eduardo Mojica-Nava}
He received a B.S. degree in Electronics Engineering from Universidad Industrial de Santander in 2002, an M.Sc. degree in Electronics and Computer Science Engineering from Universidad de Los Andes, and a Ph.D. degree in Automatique et Informatique Industrielle from École des Mines de Nantes, Nantes, France and also Universidad de Los Andes in 2010. From 2011 to 2012, he was a Post-Doctoral Researcher at Universidad de Los Andes. He has been visiting professor at Université de Mons, Mons, Belgium, and Politecnico de Milano, Milan, Italy. Currently, he is a full professor at Universidad Nacional de Colombia, Bogotá, Colombia. 
\end{IEEEbiographynophoto}







\end{document}